\documentclass{article}

\usepackage{amsmath, amssymb, amsthm, amsfonts}

\newtheorem{prop}{Proposition}
\newtheorem*{example}{Example} 

\newtheorem{theorem}{Theorem}

\newtheorem*{PL-Inequality}{Polyak-Lojasiewicz (PL) Inequality}
\DeclareMathOperator*{\argmin}{arg\,min}

\usepackage{diagbox}
\usepackage{mathtools}
\usepackage{enumerate}
\usepackage{graphicx}
\usepackage{array}
\usepackage{tikz}
\usepackage{hyperref, doi}
\usepackage{epsf}
\usepackage{bbm}
\usepackage[utf8]{inputenc} 
\usepackage[T1]{fontenc}    
\usepackage{float}
\usepackage{selectp}

\usepackage{caption}
\usepackage{subcaption}
\usepackage{array}
\usepackage{makecell}
\usepackage{authblk}
\usepackage{parskip}
\usepackage{wrapfig}

\usepackage[accepted]{icml2021}

\icmltitlerunning{Local Quadratic Convergence of Stochastic Gradient Descent with Adaptive Step Size}

\begin{document}

\twocolumn[
\icmltitle{Local Quadratic Convergence of \\ Stochastic Gradient Descent with Adaptive Step Size}
\begin{icmlauthorlist}
\icmlauthor{Adityanarayanan Radhakrishnan}{to}
\icmlauthor{Mikhail Belkin}{hc}
\icmlauthor{Caroline Uhler}{to}
\end{icmlauthorlist}
\icmlaffiliation{to}{Laboratory for Information and Decision Systems, Massachusetts Institute of Technology, Cambridge, MA, USA}
\icmlaffiliation{hc}{Halıcıoğlu Data Science Institute, University of California, San Diego, CA, USA}
\icmlcorrespondingauthor{Adityanarayanan Radhakrishnan}{aradha@mit.edu}
\icmlkeywords{Machine Learning, ICML}
\vskip 0.3in
]

\printAffiliationsAndNotice{\icmlEqualContribution} 

\begin{abstract}
Establishing a fast rate of convergence for optimization methods is crucial to their applicability in practice.  With the increasing popularity of deep learning over the past decade, stochastic gradient descent and its adaptive variants (e.g. Adagrad, Adam, etc.) have become prominent methods of choice for machine learning practitioners.  While a large number of works have demonstrated that these first order optimization methods can achieve sub-linear or linear convergence, we establish local quadratic convergence for stochastic gradient descent with adaptive step size for problems such as matrix inversion.  
\end{abstract}

\section{Introduction}

Over the past decade, gradient descent and adaptive variants such as RMSprop \citep{RmsProp} and Adam \citep{Adam} have become standard choices of optimization method due to their ability to scale to large machine learning models such as deep neural networks.  Establishing fast convergence rates and providing explicit learning rate schedules for these optimization methods has been recognized as an important problem for reducing the complexity of hyper-parameter tuning for practitioners.  Thus, analyzing the convergence of these optimization methods is an active area of research in machine learning. 

While the convergence of vanilla stochastic gradient descent has been studied extensively in the convex and non-convex setting \citep{ConvexSGD, PLConditionLinearConvergence}, understanding whether there exist adaptive step sizes that yield higher orders of convergence remains an open problem.  In this work, we show that by a careful selection of an adaptive step size, gradient descent and stochastic gradient descent can converge at a quadratic rate locally.  To build intuition, consider the following example of inverting a real number.



\begin{example}[Number Inversion]
Suppose we wish to invert a real number $x$ using only additions and multiplications.  One way to solve this problem is to find $w \in \mathbb{R}$ such that $wx = 1$ using gradient descent on the loss $\mathcal{L} =\frac{1}{2}(1- wx)^2$.  This method will yield a linear convergence rate if a fixed step size is used.  However, using an adaptive step size of ${w^{(t)}}^2$ yields a quadratic convergence rate.  The update equations from gradient descent are as follows:
\begin{align*}
    w^{(t+1)} &= w^{(t)} + \gamma^{(t)} x (1 - w^{(t)}x) \\
    &= w^{(t)} + {w^{(t)}}^2 x ( 1 - w^{(t)}x),
\end{align*}
where $\gamma^{(t)} = {w^{(t)}}^2$ is the adaptive learning rate.   Let $w^*$ be the solution to $w^*x = 1$ and let $u^{(t)} = w^{(t)} - w^*$.  Consider the convergence of $|u^{(t)}|$:
\begin{align*}
    |u^{(t+1)}| &= |u^{(t)} + (u^{(t)} + w^*)^2 x ( 1 - u^{(t)}x - w^*x)| \\
    &= |u^{(t)} - u^{(t)} (u^{(t)} + w^*)^2 x^2  | ~~~~~ \text{as $w^*x = 1$} \\
    &= |u^{(t)} (1 - {w^*}^2 x^2)  - 2{u^{(t)}}^2 w^* x^2 - {u^{(t)}}^3 x^2  | \\
    &= | - 2{u^{(t)}}^2 w^* x^2 - {u^{(t)}}^3 x^2 | \\
    &= |2{u^{(t)}}^2x + {u^{(t)}}^3 x^2 |,
\end{align*}
yielding a quadratic convergence rate for sufficiently small $u^{(0)}$.  

For this problem, Newton's method \citep[Chapter 3]{NumericalOptimizationWright} also gives a quadratic convergence rate; but interestingly, the update for Newton's method is different, namely $w^{(t+1)} = 2w^{(t)}  - {w^{(t)}}^2x$.  
\end{example}

An outline of our work is as follows.  In Section \ref{sec: Related works}, we discuss related work.  In Section \ref{sec: Adaptive GD} we provide an adaptive learning rate for local quadratic convergence of gradient descent in the multi-dimensional setting.  In Section \ref{sec: Adaptive SGD}, we extend this adaptive learning rate to establish quadratic convergence of gradient descent in the stochastic setting.  In Section \ref{sec: Experiments}, we provide experiments to demonstrate our learning rate working in practice. In Section \ref{sec: Optimality of Quadratic Convergence}, we provide a class of linear systems under which linear convergence is the highest order of convergence possible when using a learning rate that is a matrix polynomial in the parameters.  We conclude with a discussion in Section \ref{sec: Discussion}.  

\section{Related Work}
\label{sec: Related works}

Sub-linear and linear convergence rates for stochastic gradient descent (SGD) have been provided across a number of settings.  A survey of traditional results for the convex setting is presented in \citet{ConvexSGD}.  Recent works establish linear convergence of SGD in the non-convex setting \citep{BassilySGDLinearConvergence, MarkSchmidtSGDPL, PLConditionLinearConvergence} by using the Polyak-Lojasiewicz inequality \citep{PLInequalityLojasiewicz, PLInequality}.  Recent work has also focused on providing convergence rates for adaptive methods such as Adagrad \citep{Adagrad}.  In particular, \citet{SGMDRadha, AdagradPL} provide sufficient conditions for linear convergence of Adagrad and Adagrad-Norm (a normed version of Adagrad).

Several works have also developed pre-conditioner methods for gradient descent in order to improve convergence in practice, but currently, the established rates of convergence are sub-linear or linear.  For example, RMSprop \citep{RmsProp} and Adam \citep{Adam} are popular adaptive gradient methods for training deep neural networks, and these methods can be interpreted as pre-conditioner methods \citep{AdaptivePreconditioners}.  Sufficient conditions for sub-linear convergence of these methods are provided in \citet{SignSGDLearningRate, AdaptiveMethodsNonconvex}.  \citet{EigenPro} introduce EigenPro, which is a pre-conditioner for kernel machines that incorporates approximate second order information.  Recent work by \citet{SecondOrderDeepNetworks1, AdaHessian} introduce second order pre-conditioner methods for neural networks, and \citet{AdaHessian} establish a linear convergence rate for their method applied to strongly convex functions.  
 
Newton's method \citep[Chapter 3]{NumericalOptimizationWright} serves as a classical example of an algorithm that achieves quadratic convergence locally.  Recent work  presents an extension of Newton's method to compute the $p^{th}$ root (and the inverse $p^{th}$ root of a positive definite matrix) using only matrix multiplications \cite{NewtonMethodInverse}.  Under sufficient conditions on the spectrum of the matrix, \citet{SchurNewton} give an initialization under which this method converges at a quadratic rate. 

The updates in Newton's method rely on second-order information and can be expensive to compute on large datasets. Several quasi-Newton methods such as L-BFGS \citep[Chapter 3]{NumericalOptimizationWright} have been developed to avoid full construction of the Hessian, and recently, \citet{IQN} developed a stochastic quasi-Newton method with a super-linear rate of convergence.  The work by \citet{IQN} builds upon several prior works by \citet{NIMSuperlinear, SFO, stochasticQuasiNewton1, stochasticQuasiNewton2}, which also introduce stochastic quasi-Newton methods.  While the work by \citet{IQN} uses only matrix-vector multiplications like the stochastic version of our method, ours has a simpler formulation and is history-free, i.e. it does not need to store any additional terms other than the parameters.

\section{Gradient Descent with Adaptive Step Size}
\label{sec: Adaptive GD}
In the following theorem (with proof in Appendix \ref{appendix: Proof of Theorem 1}), we extend the analysis of number inversion from the introduction to the problem of inverting an (invertible) square matrix $X \in \mathbb{R}^{n \times n}$. 

\begin{theorem}
Let $X, I \in \mathbb{R}^{n \times n}$ with $I$ denoting the identity matrix, and let $X$ be invertible.  Gradient descent with right-multiplicative step size $\gamma^{(t)} = {W^{(t)}}^T W^{(t)}$ to solve $\argmin_{W} \frac{1}{2}\|I - WX \|_F^2$ converges at quadratic rate locally. 
\end{theorem}

\noindent \textbf{Intuition for the Learning Rate.} We arrived at our choice of adaptive learning rate by analyzing the gradient descent update in two-layer linear neural networks.  In particular, when using gradient descent to solve the following optimization problem
\begin{align*}
    \argmin_{W_1 \in \mathbb{R}^{d \times k}, W_2 \in \mathbb{R}^{k \times d}}\frac{1}{2}\| I - W_1W_2X\|_F^2, 
\end{align*}
the update rule for $W_1$ proceeds as follows: 
\begin{align*}
    W_1^{(t+1)} = W_1^{(t)} + (I - W_1^{(t)}W_2^{(t)}X)X^T W_2^{(t)}\gamma.
\end{align*}
Note that unlike gradient descent in the single layer setting, gradient descent in this setting involves terms of higher orders in $W_1$ and $W_2$.  Thus, we chose the learning rate appropriately to ensure that only these higher order terms remained. 
\vspace{2mm}

\noindent \textbf{Remarks.}  When the target matrix is no longer the identity but rather an invertible matrix $Y \in \mathbb{R}^{n \times n}$, we can simply transform the system such that the target is the identity.  For example, when $Y$ is an orthogonal matrix, we simply transform $X$ to $XY^T$ to reduce to the case above.  Interestingly, the update for Newton's method in this setting is given by $W^{(t+1)} = 2 W^{(t)} - W^{(t)} X W^{(t)}$, which again differs from our update rule.  Importantly, note that it is non-trivial to extend this formulation of Newton's method to the stochastic setting due to the non-commutativity of matrices.  
\vspace{2mm}

\noindent \textbf{Computing the $d^{th}$ Root of Symmetric Positive Definite Matrices.} Just as Newton's method can be extended to compute $d^{th}$ roots of matrices \citep{NewtonMethodInverse}, we can extend our method to do the same for symmetric positive definite matrices.

\begin{prop}
Let $X \in \mathbb{R}^{n \times n}$ be a symmetric positive definite matrix and let $I \in \mathbb{R}^{n \times n}$ denote the identity matrix.  Assuming $W^{(0)}$ commutes with $X$, gradient descent with right multiplicative step size $\gamma^{(t)} = \frac{1}{d^2}{W^{(t)}}^2$ to solve $\argmin_{W \in \mathbb{R}^{d}} \frac{1}{2}\|I - W^dX \|_F^2$ converges at quadratic rate locally.  
\end{prop}

The proof is presented in Appendix \ref{appendix: Proof of Proposition 1}.  Thus far, we have demonstrated that adaptive gradient descent achieves quadratic convergence for a class of problems in the full batch setting.  In the next section, we demonstrate that gradient descent with adaptive step size converges quadratically in the stochastic setting.  

\section{SGD with Adaptive Step Size}
\label{sec: Adaptive SGD}

We return to the setting of matrix inversion and now demonstrate that stochastic gradient descent with right-multiplicative learning rate ${W^{(t)}}^TW^{(t)}$ yields quadratic convergence.  
\begin{theorem}
Let $X, I\in \mathbb{R}^{n \times n}$ with $I \in \mathbb{R}^{n \times n}$ denoting the identity matrix.  Assuming all examples are seen per epoch, stochastic gradient descent with step size $\gamma^{(t)} = {W^{(t)}}^T W^{(t)}$ used to solve $\argmin_{W} \frac{1}{2}\|I - WX \|_F^2$ converges at a quadratic rate locally.  
\end{theorem}

The proof is presented in Appendix \ref{appendix: Proof of Theorem 2}.  The proof relies on the following proposition (proof in Appendix \ref{appendix: Proof of Proposition 2}), which demonstrates that using all examples per epoch eliminates the linear term in the error term, $U^{(t)}$, defined below.  

\begin{prop}
\label{prop: Prop 2} Let $W^*$ denote the solution to $WX = I$ and let $U^{(t)} = W^{(t)} - W^*$.  For positive constants $C_1, C_2$, 
$\|U^{(nt)}\|^2 \leq C_1 \| U^{(n(t-1))} \|^4 + C_2 \| U^{(n(t-1))} \|^6$. 
\end{prop}

\textbf{Remarks.} The benefit of the stochastic setting is that it involves only matrix-vector multiplications.  From the analysis above, the quadratic convergence rate in the stochastic setting is achieved by iterating through each example the same number of times per epoch.  In fact, unless all examples are iterated over the same number of times per epoch, the rate will always have a non-vanishing linear term.  While this analysis matches implementations of SGD in practice, this is different from the usual analysis of SGD where the rate is analyzed by considering the expectation over the randomness of the choice of example per step of the update \citep{MarkSchmidtSGDPL, BassilySGDLinearConvergence}.  Proceeding according to this usual style of analysis will only lead to an expected linear convergence rate.  

\section{Experiments}
\label{sec: Experiments}
In the following, we demonstrate that our gradient descent and stochastic gradient descent algorithms also achieve the quadratic rate in practice.  We begin by showing that when initialized sufficiently close to the true solution $W^*$ (i.e., satisfying the condition on $U^{(0)}$), our algorithms achieve quadratic convergence when inverting large random matrices.  

\par{\textbf{Construction of Invertible Matrix $X$.}} The singular values of $X$ are created by sampling $n$ values from a standard normal distribution and taking absolute value. Then the left and right orthonormal matrices $U, V$ of $X$ are generated by taking $U, V$ from the SVD of a random matrix with i.i.d.~standard normal entries.   

\par
\begin{figure*}[!t]
    \centering
    \includegraphics[height=2in]{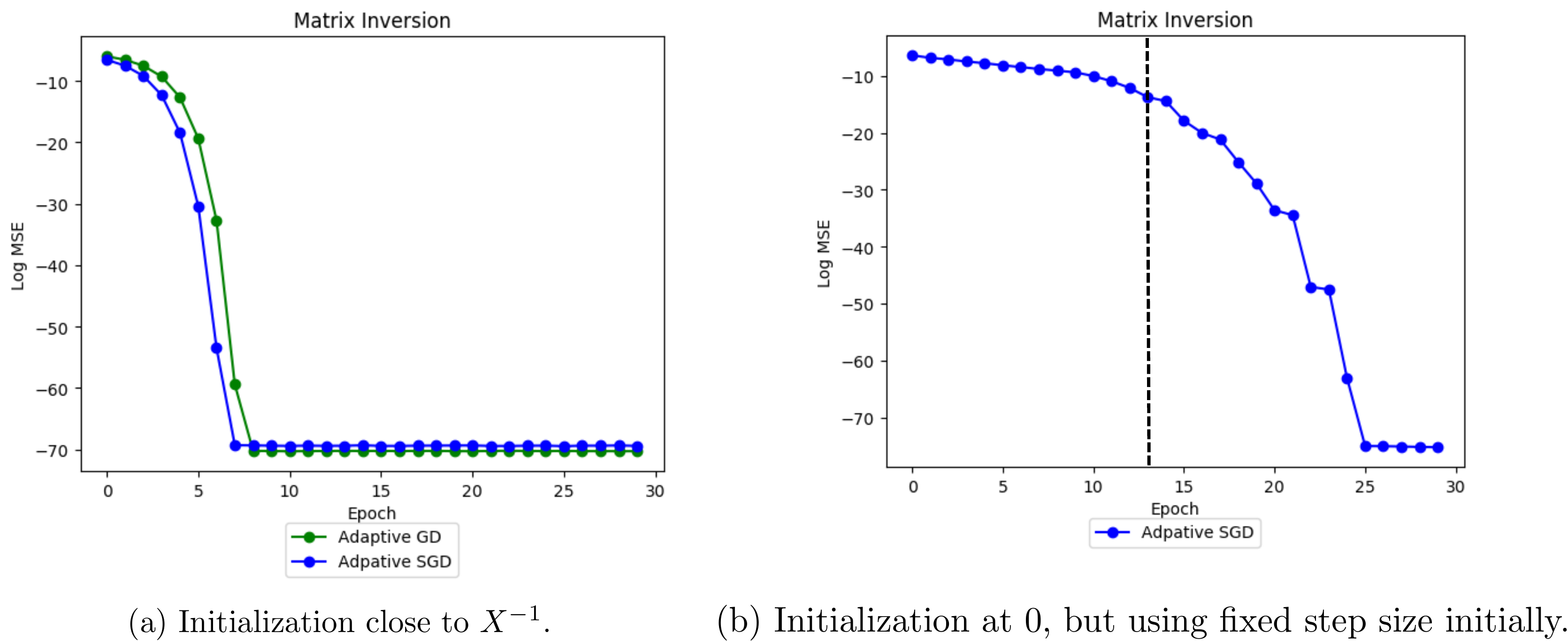}
    \caption{Using gradient descent with adaptive step size leads to locally quadratic convergence rate. In (a), we demonstrate quadratic convergence when initialization is close to the true solution.  In (b), we demonstrate that for well conditioned matrices, standard gradient descent (with linear rate) can be run until the learned solution is close enough to $X^{-1}$ (indicated by the vertical, dashed line), at which point the quadratic rate can be substituted to accelerate convergence.}
    \label{fig:Matrix Inversion 100x100}
\end{figure*}

{\textbf{Quadratic Convergence of GD with Adaptive Step Size.}}  In Figure \ref{fig:Matrix Inversion 100x100}a, we demonstrate our algorithms on inverting a random $100 \times 100$ matrix $X$.  In Figure \ref{fig:Matrix Inversion 100x100}a, we have initialized $W^{(0)}$ for gradient descent and stochastic gradient descent to be $.4 X^{-1}$ and $.5 X^{-1}$ respectively\footnote{Initializing both with $.5 X^{-1}$ yields overlapping curves.}.  When initializing close to the true solution $X^{-1}$, we see that both algorithms achieve quadratic convergence rates, in line with our theoretical results.

Of course in practice, one cannot initialize in this manner as it requires computing $X^{-1}$. We demonstrate in Figure \ref{fig:Matrix Inversion 100x100}b that we can accelerate matrix inversion by using usual gradient descent (or alternatively randomized Kaczmarz iteration \cite{KaczmarzLinearConvergence}) with initialization at $0$ until the loss is sufficiently low, and then substitute our adaptive algorithm to accelerate convergence.  In Figure \ref{fig:Matrix Inversion 100x100}b, we invert $X$ by first using gradient descent with learning rate $0.1$ and then after the loss is below $10^{-4}$, we switch to our adaptive method, which converges quadratically.    

\par
\begin{figure*}[!t]
    \centering
    \includegraphics[height=2in]{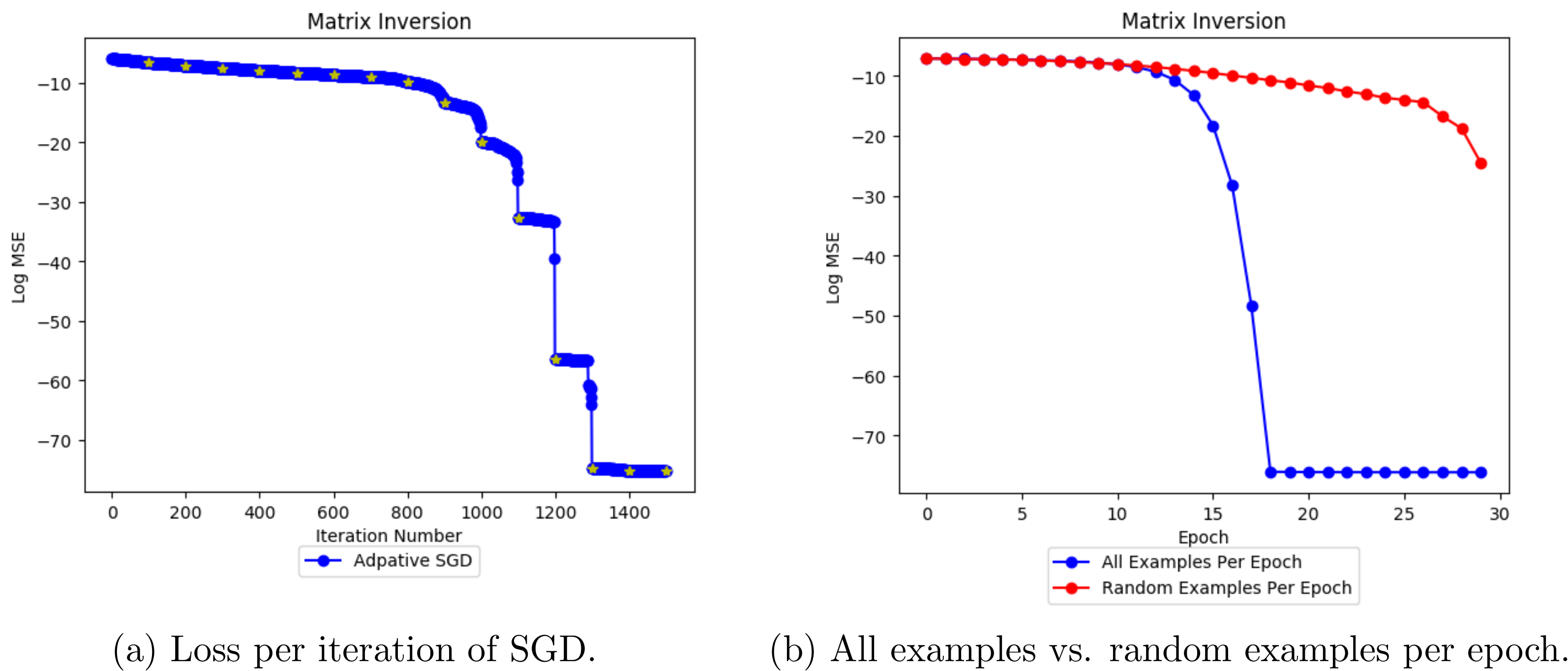}
    \caption{Stochastic gradient descent with adaptive step size converges quadratically if and only if an epoch consists of an equal count of each of the training examples. In (a), we demonstrate that within a single epoch our adaptive SGD algorithm has linear convergence rate until the last few examples are seen.  In (b), we demonstrate that using randomly sampled data instead of using all examples per epoch can result in linear convergence.}
    \label{fig:Matrix Inversion 100x100 SGD}
\end{figure*}

{\textbf{Quadratic Convergence of SGD with Adaptive Step Size.}} We now demonstrate that for our Adaptive SGD method, each of the $n$ examples must be used the same number of times per epoch in order to achieve quadratic convergence.  Instead of visualizing the log-loss per epoch, we visualize the loss per iteration in Figure \ref{fig:Matrix Inversion 100x100 SGD}a.  In this figure the stars indicate the loss after an epoch.  There are sharp drops in the loss function due to individual examples, as seen in the example used near iteration 1200.  Without the inclusion of this example, the loss would decrease linearly instead of quadratically.  In fact, if we randomly sample examples at each iteration instead of cycling through all examples per epoch, we can only guarantee a linear convergence rate.  This is demonstrated in Figure \ref{fig:Matrix Inversion 100x100 SGD}b in which we invert a random $1000 \times 1000$ matrix using adaptive SGD with randomly selected samples per epoch versus using all examples per epoch.  Both methods are initialized close to the true solution (at $.1 X^{-1}$), but randomly sampling per epoch can lead to a linear rate of convergence.

\section{What is the Optimal Rate of Convergence?}
\label{sec: Optimality of Quadratic Convergence}

While the above examples demonstrate that gradient descent with adaptive step size can converge at a quadratic rate locally, is there an alternate learning rate schedule under which we can establish an even faster rate of convergence? In the following theorem (proof in Appendix \ref{appendix: Proof of Theorem 3}), we prove that for linear systems, gradient descent with matrix polynomial learning can achieve at most linear convergence when the solution $W^* \in \mathbb{R}^{k \times d}$ has rank less than~$d$.  


\begin{theorem}
\label{thm: Theorem 3}
Let $X \in \mathbb{R}^{d \times n}$ be generic, $Y \in \mathbb{R}^{k \times n}$, and suppose that there exists $W^* \in \mathbb{R}^{k \times d}$ such that $W^* X = Y$ and $W^*$ has rank less than $d$.  Gradient descent with matrix polynomial learning rate $\gamma^{(t)} = P[{W^{(t)}}^T W^{(t)}]$ used to solve $\argmin_{W \in \mathbb{R}^{k\times d}} \frac{1}{2} \|Y - W X\|_F^2$ converges with order at most 1.
\end{theorem}

\vspace{-3mm}
\textbf{Remarks.} We note that Theorem \ref{thm: Theorem 3} does not rule out the possibility that there may be adaptive learning rates involving polynomials of sufficiently high degree in $X, Y,$ and $W^{(t)}$ that can provide higher order convergence.  However, the cost of computing such polynomials must also be taken into account when considering the practical application of such a method.  For example, if the $c_i$ terms in the proof of Theorem \ref{thm: Theorem 3} are matrices, and $\gamma^{(t)}$ is a matrix polynomial of sufficiently high degree $r$, then we can apply linear regression to find the coefficients $c_i$ to provide as fast a convergence rate as desired.  However, solving a regression problem to find these coefficients is impractical.  

\section{Discussion}
\label{sec: Discussion}
In this work, we first demonstrated that gradient descent with adaptive step size can converge at a quadratic rate locally.  We then extended our result to the stochastic setting and verified our results empirically.  Lastly, we demonstrated that if the learning rate is a matrix polynomial in the parameters and the number of regression targets is fewer than the number of features, the best possible rate of convergence is linear.  We now briefly discuss future directions of interest.  

While our quadratic convergence results are local, it would be interesting to understand whether there is a simple initialization scheme (as in \citet{SchurNewton}) under which our method achieves this rate.  Another avenue of interest would be to understand whether there are other computationally tractable learning rates that could provide even higher orders of convergence locally.  

\section*{Acknowledgements}
A.~Radhakrishnan and C.~Uhler were partially supported by NSF (DMS-1651995), ONR (N00014-17-1-2147 and N00014-18-1-2765), the MIT-IBM Watson AI Lab, and a Simons Investigator Award to C. Uhler. The Titan Xp used for this research was donated by the NVIDIA Corporation. M.~Belkin acknowledges support from NSF (IIS-1815697 and IIS-1631460) and a Google Faculty Research Award.  
\bibliographystyle{icml2021}
\bibliography{references}

\begin{thebibliography}{26}
\providecommand{\natexlab}[1]{#1}
\providecommand{\url}[1]{\texttt{#1}}
\expandafter\ifx\csname urlstyle\endcsname\relax
  \providecommand{\doi}[1]{doi: #1}\else
  \providecommand{\doi}{doi: \begingroup \urlstyle{rm}\Url}\fi

\bibitem[Anil et~al.(2020)Anil, Gupta, Koren, Regan, and
  Singer]{SecondOrderDeepNetworks1}
Anil, R., Gupta, V., Koren, T., Regan, K., and Singer, Y.
\newblock Second order optimization made practical.
\newblock \emph{arXiv preprint arXiv:2002.09018}, 2020.

\bibitem[Bassily et~al.(2018)Bassily, Belkin, and
  Ma]{BassilySGDLinearConvergence}
Bassily, R., Belkin, M., and Ma, S.
\newblock On exponential convergence of {SGD} in non-convex over-parametrized
  learning.
\newblock \emph{arXiv preprint arXiv:1811.02564}, 2018.

\bibitem[Bernstein et~al.(2018)Bernstein, Wang, Azizzade-nesheli, and
  Anandkumar]{SignSGDLearningRate}
Bernstein, J., Wang, Y.-X., Azizzade-nesheli, K., and Anandkumar, A.
\newblock {signSGD: Compressed optimisation for non-convex problems}.
\newblock In \emph{International Conference on Machine Learning}, 2018.

\bibitem[Bubeck(2015)]{ConvexSGD}
Bubeck, S.
\newblock {Convex optimization: Algorithms and complexity}.
\newblock \emph{Foundations and Trends in Machine Learning}, 8\penalty0
  (3-4):\penalty0 231--357, 2015.

\bibitem[Byrd et~al.(2016)Byrd, Hansen, Nocedal, and
  Singer]{stochasticQuasiNewton2}
Byrd, R., Hansen, S., Nocedal, J., and Singer, Y.
\newblock Stochastic {quasi-Newton} methods for nonconvex stochastic
  optimization.
\newblock \emph{SIAM Journal on Optimization}, 26\penalty0 (2):\penalty0
  1008--1031, 2016.

\bibitem[Duchi et~al.(2011)Duchi, Hazan, and Singer]{Adagrad}
Duchi, J., Hazan, E., and Singer, Y.
\newblock {Adaptive subgradient methods for online learning and stochastic
  optimization}.
\newblock \emph{Journal of Machine Learning Research}, 12:\penalty0 2121--2159,
  2011.

\bibitem[Guo \& Higham(2006)Guo and Higham]{SchurNewton}
Guo, C.-H. and Higham, N.~J.
\newblock {A Schur-Newton method for the matrix $p$th root and its inverse.}
\newblock \emph{SIAM Journal on Matrix Analysis and Applications}, 28\penalty0
  (3):\penalty0 788–804, 2006.

\bibitem[Iannazzo(2005)]{NewtonMethodInverse}
Iannazzo, B.
\newblock {On the Newton method for the matrix $P$th root}.
\newblock \emph{SIAM Journal on Matrix Analysis and Applications}, 28\penalty0
  (2):\penalty0 503--523, 2005.

\bibitem[Karimi et~al.(2016)Karimi, Nutini, and
  Schmidt]{PLConditionLinearConvergence}
Karimi, H., Nutini, J., and Schmidt, M.
\newblock Linear convergence of gradient and proximal-gradient methods under
  the {Polyak-Lojasiewicz} condition.
\newblock In \emph{Joint European Conference on Machine Learning and Knowledge
  Discovery in Databases}, pp.\  795--811. Springer, 2016.

\bibitem[Kingma \& Ba(2015)Kingma and Ba]{Adam}
Kingma, D.~P. and Ba, J.
\newblock {Adam: A method for stochastic optimization}.
\newblock In \emph{International Conference on Learning Representations}, 2015.

\bibitem[Lojasiewicz(1963)]{PLInequalityLojasiewicz}
Lojasiewicz, S.
\newblock {A topological property of real analytic subsets (in French)}.
\newblock \emph{Les \'equations aux d\'eriv\'ees partielles.}, 117:\penalty0
  87--89, 1963.

\bibitem[Ma \& Belkin(2017)Ma and Belkin]{EigenPro}
Ma, S. and Belkin, M.
\newblock {Diving into the shallows: a computational perspective on large-scale
  shallow learning}.
\newblock In \emph{Advances in Neural Information Processing Systems}, 2017.

\bibitem[Mokhtari et~al.(2018)Mokhtari, Eisen, and Ribeiro]{IQN}
Mokhtari, A., Eisen, M., and Ribeiro, A.
\newblock {IQN: An incremental quasi-Newton method with local superlinear
  convergence rate}.
\newblock \emph{SIAM Journal on Optimization}, 28\penalty0 (2):\penalty0
  1670--1698, 2018.

\bibitem[Nocedal \& Wright(2006)Nocedal and
  Wright]{NumericalOptimizationWright}
Nocedal, J. and Wright, S.~J.
\newblock \emph{{Numerical Optimization}}, volume~2.
\newblock Springer, 2006.

\bibitem[Polyak(1963)]{PLInequality}
Polyak, B.
\newblock {Gradient methods for minimizing functionals (in Russian)}.
\newblock \emph{Zh. Vychisl. Mat. Mat. Fiz.}, 3\penalty0 (4):\penalty0
  643–653, 1963.

\bibitem[Radhakrishnan et~al.(2020)Radhakrishnan, Belkin, and Uhler]{SGMDRadha}
Radhakrishnan, A., Belkin, M., and Uhler, C.
\newblock Linear convergence and implicit regularization of generalized mirror
  descent with time-dependent mirrors.
\newblock \emph{arXiv preprint arXiv:2009.08574}, 2020.

\bibitem[Rodomanov \& Proezd(2016)Rodomanov and Proezd]{NIMSuperlinear}
Rodomanov, A. and Proezd, K.
\newblock A superlinearly-convergent proximal {Newton-type} method for the
  optimization of finite sums.
\newblock In \emph{International Conference on Machine Learning}, 2016.

\bibitem[Sohl-Dickstein et~al.(2014)Sohl-Dickstein, Poole, and Ganguli]{SFO}
Sohl-Dickstein, J., Poole, B., and Ganguli, S.
\newblock {Fast large-scale optimization by unifying stochastic gradient and
  {quasi-Newton} methods}.
\newblock In \emph{International Conference on Machine Learning}, 2014.

\bibitem[Staib et~al.(2020)Staib, Reddi, Kale, Kumar, and
  Sra]{AdaptivePreconditioners}
Staib, M., Reddi, S., Kale, S., Kumar, S., and Sra, S.
\newblock Escaping saddle points with adaptive gradient methods.
\newblock \emph{arXiv preprint arXiv:1901.09149}, 2020.

\bibitem[Strohmer \& Vershynin(2009)Strohmer and
  Vershynin]{KaczmarzLinearConvergence}
Strohmer, T. and Vershynin, R.
\newblock A randomized {Kaczmarz} algorithm with exponential convergence.
\newblock \emph{Journal of Fourier Analysis and Applications}, 15\penalty0
  (292), 2009.

\bibitem[Tieleman \& Hinton(2012)Tieleman and Hinton]{RmsProp}
Tieleman, T. and Hinton, G.
\newblock {Lecture 6.5---RmsProp: Divide the gradient by a running average of
  its recent magnitude}.
\newblock COURSERA: Neural Networks for Machine Learning, 2012.

\bibitem[Vaswani et~al.(2019)Vaswani, Bach, and Schmidt]{MarkSchmidtSGDPL}
Vaswani, S., Bach, F., and Schmidt, M.
\newblock Fast and faster convergence of {SGD} for over-parameterized models
  and an accelerated perceptron.
\newblock In \emph{International Conference on Artificial Intelligence and
  Statistics}, 2019.

\bibitem[Wang et~al.(2017)Wang, Ma, Goldfarb, and Liu]{stochasticQuasiNewton1}
Wang, X., Ma, S., Goldfarb, D., and Liu, W.
\newblock Stochastic {quasi-Newton} methods for nonconvex stochastic
  optimization.
\newblock \emph{SIAM Journal on Optimization}, 27\penalty0 (2):\penalty0
  927--956, 2017.

\bibitem[Xie et~al.(2020)Xie, Wu, and Ward]{AdagradPL}
Xie, Y., Wu, X., and Ward, R.
\newblock Linear convergence of adaptive stochastic gradient descent.
\newblock In \emph{International Conference on Artificial Intelligence and
  Statistics}, 2020.

\bibitem[Yao et~al.(2020)Yao, Gholami, Shen, Keutzer, and Mahoney]{AdaHessian}
Yao, Z., Gholami, A., Shen, S., Keutzer, K., and Mahoney, M.~W.
\newblock {ADAHESSIAN: An} adaptive second order optimizer for machine
  learning.
\newblock \emph{arXiv preprint arXiv:2006.00719}, 2020.

\bibitem[Zaheer et~al.(2018)Zaheer, Reddi, Sachan, Kale, and
  Kumar]{AdaptiveMethodsNonconvex}
Zaheer, M., Reddi, S., Sachan, D., Kale, S., and Kumar, S.
\newblock {Adaptive methods for nonconvex optimization}.
\newblock In \emph{Advances in Neural Information Processing Systems}, 2018.

\end{thebibliography}

\appendix
\renewcommand{\thesubsection}{\Alph{subsection}}
\onecolumn

\section{Proof of Theorem 1}
\label{appendix: Proof of Theorem 1}
\begin{proof}
Gradient descent proceeds as follows:
\begin{align*}
    W^{(t+1)} &= W^{(t)} +  (I - W^{(t)}X) X^T \gamma^{(t)} \\
    &= W^{(t)} +  (I - W^{(t)}X) X^T {W^{(t)}}^T W^{(t)}.
\end{align*}
Note that $\gamma^{(t)}$ needs to be multiplied on the right hand side of the gradient, i.e., it is a non-commutative learning rate.  Let $U^{(t)} = W^{(t)} - W^*$ where $W^*$ is the solution satisfying $W^* X = I$, and let $\sigma_1 \geq \sigma_2 \ldots \geq \sigma_n$ be the singular values of $X$.  To simplify notation, we assume that the norms used are Frobenius norms unless otherwise specified.  Then we obtain:
\begin{align*}
    \| U^{(t+1)} \|^2 &= \| U^{(t)} +  ( I - U^{(t)}X - W^*X)X^T (U^{(t)} + W^*)^T (U^{(t)} + W^*) \|^2\\
    &= \| U^{(t)} - U^{(t)} X X^T ({U^{(t)}}^T U^{(t)} + {W^*}^T U^{(t)} +  {U^{(t)}}^T W^* + {W^*}^T W^*) \|^2 ~ \text{as $W^*X = I$}  \\
    &= \| U^{(t)} (I -  X X ^T {W^*}^T W^*)  - U^{(t)} X X^T ({U^{(t)}}^T U^{(t)} + {W^*}^T U^{(t)} +  {U^{(t)}}^T W^*) \|^2  \\
    &= \|- U^{(t)} X X^T ({U^{(t)}}^T U^{(t)} + {W^*}^T U^{(t)} +  {U^{(t)}}^T W^*) \|^2 \\
    &\leq \|X X^T \|_2^2 (({\| U^{(t)}\|^2})^3  + \|W^* \|^2 {\|U^{(t)}\|^2}^2) \\
    &\leq \sigma_1^4 {\| U^{(t)}\|^6}  + n \sigma_1^4 \sigma_n^{-2} {\|U^{(t)}\|^4}.
\end{align*}
Hence provided that $\| U^{(0)} \|^2$ is sufficiently close to zero, $\| U^{(t)} \|^2$ converges quadratically.  
\end{proof}

\section{Proof of Proposition 1}
\label{appendix: Proof of Proposition 1}
\begin{proof}
By assumption, $W^{(0)}$ commutes with $X$ so that $W^{(t)}$ is a matrix polynomial in $X$ and hence commutes with $X$ for all $t$.  In particular, the update equation can be written as
\begin{align*}
    W^{(t+1)} &= W^{(t)} + d \gamma^{(t)} {W^{(t)}}^{d-1}(I - {W^{(t)}}^d X) X \\
    &= W^{(t)} + \frac{1}{d} {W^{(t)}}^{d+1}(I - {W^{(t)}}^d X) X.
\end{align*}   
We let $U^{(t)} = W^{(t)} - W^*$ where ${W^*}^d X = I$.  As before, we obtain 
\begin{align*}
    \|U^{(t+1)}\|^2 &= \|U^{(t)} + \frac{1}{d}(U^{(t)} + W^*)^{d+1} (I - (U^{(t)} + W^*)^d X )X \|^2 \\
    &\leq  \| U^{(t)}[I - {W^*}^{d+1} {W^*}^{d-1} X^2] \|^2 + \displaystyle\sum\limits_{i=2}^{2d+1} C_i \| {U^{(t)}}^i \|^2  ~~~~~ \text{for constants $C_i \in \mathbb{R}_+$}\\
    &= \displaystyle\sum\limits_{i=2}^{2d+1} C_i \| {U^{(t)}}^i \|^2 f\leq \displaystyle\sum\limits_{i=2}^{2d+1} C_i \| {U^{(t)}} \|^{2i}.
\end{align*}

Hence provided that $\| U^{(0)} \|^2$ is sufficiently small,  gradient descent with adaptive step size achieves a quadratic convergence rate.  
\end{proof}

\section{Proof of Theorem 2}
\label{appendix: Proof of Theorem 2}
\begin{proof}
Let $X_i$ denote the $i^{th}$ column vector of $X$ and let $e_i$ denote the $i^{th}$ column vector of $I$ (that is $e_i$ is a vector with a $1$ in the $i^{th}$ coordinate and $0$ elsewhere).  We let an epoch be a pass through a random ordering of the $n$ pairs $(X_i, e_i)$ where each pair corresponds to the update
\begin{align}
\label{eqn: SGD update}
    W^{(t+1)} = W^{(t)} + (e_i - W^{(t)}X_i) X_i^T \gamma^{(t)},
\end{align}
where we set $\gamma^{(t)} = {W^{(t)}}^T W^{(t)}$. Defining $U^{(t)} = W^{(t)} - W^*$, we will show that the sequence $\| U^{(nt)} \|^2$ for $t \in \mathbb{N}$ converges at a quadratic rate.  We prove this by using Proposition \ref{prop: Prop 2} below, which states that the term of $U^{(nt)}$ that is linear in $U^{n(t-1)}$ has coefficient $0$.  In particular, Proposition \ref{prop: Prop 2} implies that $\| U^{(nt)} \|^2 \leq C_1 {\| U^{(n(t-1))} \|^4} + C_2{\| U^{(n(t-1))} \|^6}$ for constsants $C_1, C_2 > 0$ .  Hence, initializing such that  $\|U^{(0)}\|^2$ is sufficiently close to $0$ yields a quadratic rate of convergence.   
\end{proof}

\section{Proof of Proposition 2}
\label{appendix: Proof of Proposition 2}

\begin{proof}
We prove this for the case when $t = 1$, from which the general result follows easily.  Without loss of generality, let the order of updates for the first epoch be $(X_1, e_1), (X_2, e_2) \ldots (X_n, e_n)$.  We first show that $C_0 = \| \displaystyle\prod\limits_{i=1}^{n} (I - X_i X_i^T {W^*}^T W^*)  \|^2$.  

\noindent By equation \ref{eqn: SGD update}, we have that
\begin{align*}
    U^{(i)} =&\; U^{(i-1)} + (e_i - (U^{(i-1)} + W^*)X_i)  X_i^T (U^{(i-1)} + W^*)^T (U^{(i-1)} + W^*) \\
    =&\; U^{(i-1)} - U^{(i-1)} X_i X_i^T ({U^{(i-1)}}^T U^{(i-1)} + {U^{(i-1)}}^T W^* + {W^*}^T U^{(i-1)}  + {W^*}^T W^*) \\ 
    =&\; U^{(i-1)} [I - X_i X_i^T {W^*}^T W^*]  - U^{(i-1)} X_i X_i^T {U^{(i-1)}}^T W^* - U^{(i-1)} X_i X_i^T {W^*}^T U^{(i-1)} \\&\; - U^{(i-1)} X_i X_i^T {U^{(i-1)}}^T U^{(i-1)}. 
\end{align*}
Hence the coefficient of $U^{(i-1)}$ in the expansion of $U^{(i)}$ is $[I - X_i X_i^T {W^*}^T W^*]$.  Proceeding recursively from $i = n$ to $i=1$ proves that the coefficient of $U^{(0)}$ in the expansion of $U^{(n)}$ is $\displaystyle\prod\limits_{i=1}^{n} (I - X_i X_i^T {W^*}^T W^*)$.

\noindent Now we show that this term is $\mathbf{0}$. First note that $W^* X_i = e_i$ and $X_i^T {W^*}^T = {e_i}^T$.  Hence, we obtain that for $i \neq j$,
\begin{align*}
    X_i X_i^T {W^*}^T W^* X_j X_j^T {W^*}^T W^* &= X_i (X_i^T {W^*}^T) (W^* X_j) X_j^T {W^*}^T W^* \\
    &= X_i e_i^T e_j X_j^T {W^*}^T W^* \\
    &= \mathbf{0} ~~~~ \text{as $e_i^T e_j = 0$ for $i \neq j$},
\end{align*}
and thus that
\begin{align*}
    \displaystyle\prod\limits_{i=1}^{n} (I - X_i X_i^T {W^*}^T W^*) &= I - \displaystyle\sum\limits_{i=1}^{n} X_i X_i^T {W^*}^T W^* 
    = I - \displaystyle\sum\limits_{i=1}^{n} X_i e_i^T W^* 
    = I - X W^*
    = \mathbf{0}, 
\end{align*}
which completes the proof.
\end{proof}

\section{Proof of Theorem 3}
\label{appendix: Proof of Theorem 3}
\begin{proof}
    Gradient descent proceeds according to the update:
    \begin{align*}
        W^{(t+1)} &= W^{(t)} +  (Y - W^{(t)} X) X^T \gamma^{(t)}.
    \end{align*}
    Again, assuming that $U^{(t)} = W^{(t)} - W^*$, we have:
    \begin{align*}
        U^{(t+1)} &= U^{(t)} -  U^{(t)} XX^T \gamma^{(t)} = U^{(t)}[ I_{d} - XX^T \gamma^{(t)}].
    \end{align*}
    Now letting $\gamma^{(t)} = \sum_{i=0}^{r} c_i ({W^{(t)}}^T W^{(t)})^i $, we obtain
    \begin{align*}
        \gamma^{(t)} = \sum_{i=0}^{r} c_i ((U^{(t)} + W^*)^T (U^{(t)} + W^*))^i.
    \end{align*}
Note that the constant term of $U^{(t)}$ in $\gamma^{(t)}$ is given by $C_0 = \sum_{i=0}^{r}c_i ({W^*}^T W^*)^i$ for $\{c_i\}_{i=1}^{r} \subset \mathbb{R}$.  Thus, if $W^*$ has rank less than $d$, then $XX^T C_0$ cannot equal $I_d$.  Hence the term of $U^{(t+1)}$ that is linear in $U^{(t)}$ has a non-zero coefficient, which completes the proof. 
\end{proof}

\end{document}